\documentclass[12pt]{amsart}
\usepackage{amsfonts,amssymb,hyperref}
\newtheorem{thm}{Theorem}[section]

\newtheorem{lem}[thm]{Lemma}
\newtheorem{prop}[thm]{Proposition}
\theoremstyle{definition}

\theoremstyle{remark}
\newtheorem{rem}[thm]{Remark}
\numberwithin{equation}{section}
\newcommand{\norm}[1]{\left\Vert#1\right\Vert}
\newcommand{\abs}[1]{\left\vert#1\right\vert}

\newcommand{\Real}{\mathbb R}

\newcommand{\To}{\longrightarrow}

\newcommand{\Ind}{{\rm Ind}}
\newcommand\Det{\mathop{\rm Det}}

\def\m{\underline{\bf m}}
\def\n{\underline{\bf n}}
\begin{document}

\title[Ring structures for discrete series]{Ring structures for holomorphic discrete series and Rankin-Cohen brackets}
\author{Gerrit van DIJK, Michael PEVZNER}%
\address{G. van Dijk: Mathematish Instituut, Universiteit Leiden, PO Box 9512, NL-2300 RA
Leiden, Nederland.}\email{dijk@math.leidenuniv.nl}
\address{
M.Pevzner: Laboratoire de Math\'ematiques, UMR CNRS 6056, Universit\'e de Reims, Campus Moulin de la Housse BP 1039, F-51687, Reims, France.}%
\email{pevzner@univ-reims.fr}%


\begin{abstract}
In the present note we discuss two different ring structures on the
set of holomorphic discrete series of a causal symmetric space of
Cayley type $G/H$ and we suggest a new interpretation of
Rankin-Cohen brackets in terms of intertwining operators arising in
the decomposition of tensor products of holomorphic discrete series representations.
\end{abstract}
\maketitle

\section{Introduction}
When studying $L$-functions of quadratic characters H. Cohen
\cite{[Cohen]} described in 1975 a particular family of
bi-differential operators acting on smooth functions on the
Poincar\'e upper half-plane $\Pi$. The initial interest in these
operators, called henceforth the Rankin-Cohen brackets (RCB), is due
to the fact that they give a powerful tool for producing new modular
forms of higher weight.

More precisely, fix a positive integer $k$ and define for every
$f\in C^\infty(\Pi)$ :
$$
(f_{|k}\gamma)(z):=(cz+d)^{-k}f\left(\frac{az+b}{cz+d}\right),\qquad
\forall\gamma=\left(%
\begin{array}{cc}
  a & b \\
  c & d \\
\end{array}%
\right)\in SL(2,\Real).
$$
One says that a function $f$ holomorphic on $\Pi$ is a modular form
of weight $k$  with respect to some arithmetic subgroup $\Gamma\subset
G$ if it satisfies the identity $(f_{|k}\gamma)=f$ for all
$\gamma\in\Gamma$.

Let $k_1,k_2,j$ be three positive integers and $f,g\in
C^\infty(\Pi)$. One sets
\begin{equation}\label{rcb}
F_j(f,g)=\sum_{\ell=0}^j(-1)^{\ell}C_{k_1+j-1}^{\ell}C_{k_2+j-1}^{j-\ell}f^{(j-\ell)}g^{(\ell)},\:
{\mathrm {where}}\: f^{(\ell)}=\left(\frac{\partial}{\partial
z}\right)^{\ell}f,
\end{equation}
and $C_k^\ell$ denote the binomial coefficient $\frac{k!}{(k-\ell)!\ell!}$.

H. Cohen showed that the following identity holds:
\begin{equation*}\label{identite}
    F_j(f_{|k_1}\gamma,g_{|k_2}\gamma)=F_j(f,g)_{|k_1+k_2+2j}\gamma,\qquad\gamma\in
    SL(2,\mathbb R).
\end{equation*}
Therefore if $f$ and $g$ are modular of weight $k_1$ and $k_2$
respectively $F_j(f,g)$ is again a modular form of weight
$k_1+k_2+2j$ for every $j\in\mathbb N$. Notice that in case when
$\Gamma=SL(2,\mathbb Z)$ the only non trivial modular forms are of
even weight.\vskip10pt

This construction was generalized in the setting of $Sp(n,\mathbb
Z)$-modular forms on the Siegel half plane, i.e. the symmetric space
of positive definite symmetric matrices, by W.Eholzer and T.Ibukiyama
\cite{[Eholzer]}.
 An algebraic approach to RCB and their possible
generalizations via the commutation relations that they should
satisfy was developed by D. Zagier \cite{[Zagier1]}. See also \cite{[Zagier2]}
for an overview of this subject from the number theoretic point of view.
\bigskip

From the other side, in 1996 A and J. Unterberger \cite{[UU]} showed
that this family of bi-differential operators arises in an
astonishing way in the context of the covariant quantization of one
of the coadjoint orbits of the Lie group $G=SL(2,\Real)$. By developing
a covariant symbolic calculus on the one-sheeted hyperboloid
realized as the symmetric space $G/H=SL(2,\Real)/SO(1,1)$ they proved
that the composition $f\#_s g$ of two symbols $f$ and $g$
satisfying some regularity conditions  (they are images by the
inverse Laplace transform of holomorphic functions, square
integrable with respect to some particular measure on the upper half plane)
is again a symbol of the same kind and moreover it decomposes into a
convergent sum $f\#_s g=\sum_jh_j$ where every summand $h_j$ is
related to the Rankin-Cohen bracket $F_j(f,g)$.

This result implies that the set of holomorphic discrete series
representations with even parameter of the group $SL(2,\Real)$ is endowed with
a graded non-commutative ring structure given by the so-called
standard (or convolution-first) covariant symbolic calculus on
$SL(2,\Real)/SO(1,1)$.

\bigskip

The group of unimodular real matrices $G=SL(2,\mathbb{R})$ acts on
the set of functions defined on $\Pi$ and the modular forms are the
invariants of this action restricted to $SL(2,\mathbb Z)\subset G$.

The fact that RCB's produce new modular forms from known ones fits
with the standard techniques of transvectants developed in the
classical invariant theory. This method allows us to construct new
invariant analytic functions in two complex variables starting with
a couple of known analytic functions invariant for the simultaneous
linear action of $GL(2,\mathbb C)$. This procedure involves some
differential operators such that once restricted to homogeneous
functions they coincide with the RCB's given by (\ref{rcb}). P.
Olver gives a very detailed overview of this construction in chapter
5 of his book \cite{[Olver]} as well as in \cite{[OlverSan]}. Notice
that the basic lemma underlying the link between transvectants ant
the Rankin-Cohen brackets was proved by S. Gundelfinger
\cite{[Gund]} already in 1886.\vskip10pt

Inspired by this observation we shall gather in the present note
these two different approaches to the RCB's using the representation
theory of the group $SL(2,\mathbb R)$. These techniques will make
clear the way to generalize the notion of RCB's in the setting of
para-Hermitian symmetric spaces of Hermitian type (some times called
also causal symmetric spaces of Cayley type). The choice of this
particular class of symplectic symmetric spaces is explained in the
next section. We shall see that RCB's are related to the
decomposition of tensor products of two holomorphic discrete series
representations into irreducible components. Recent results by Peetre \cite{[Peetre]} and Peng and Zhang
\cite{[PZ]} give an explicit formula for the RCB in this case. The description of the
Clebsh-Gordan coefficients of the group $G$ is an important problem
even from the physical point of view and we hope that this note will
give a better understanding of what one calls now the Rankin-Cohen
quantization \cite{[CM]}.
\bigskip

 \thanks{M.P. is grateful to J. Alev, A. Unterberger and G. Zhang for fruitful
 discussions and thanks the E. Schr\"odinger International Institute
 for Mathematical Physics in Vienna for its hospitality and support.}

\section{Geometric settings}

Let $G$ be a connected real semi-simple Lie group with finite center
and $K$ be a maximal compact subgroup. We assume that the
Harish-Chandra condition ${\rm rank}\, G={\rm rank}\, K$ holds what
guarantees the existence of discrete series representations (i.e.
unitary irreducible representations whose matrix coefficients are
square integrable on $G$). Furthermore, we assume that $G/K$ is a
Hermitian symmetric space of tube type and thus $G$ has holomorphic
discrete series, i.e. discrete series realizable in holomorphic
sections of holomorphic vector bundles over $G/K$. Equivalently the
last condition means that the Harish-Chandra modules underlying these
discrete series representations are highest weight modules.

Among such Lie groups we shall restrict our considerations to those
which can be seen as automorphism groups of some semi-simple
para-Hermitian symmetric space. More precisely, let $\sigma$ be an
involutive automorphism of $G$ and $H$ an open connected subgroup of
the group of fixed points of $\sigma$. The coset space $G/H$ (which
is actually a coadjoint orbit of $G$ and therefore is a symplectic
manifold) is called para-Hermitian if its tangent bundle $T(G/H)$
splits into the sum of two $G$-invariant isomorphic sub-bundles (see
\cite{[KK]} for a detailed study of such spaces).

This splitting induces a $G$-invariant polarization on $T(G/H)$
which is necessary in order to define a symbolic calculus. Indeed,
this polarization will allow us to distinguish position and momenta
variables on $G/H$ which plays the role of the phase space, while
the symmetric space $G/K$ will be seen as the configuration space.

It turns out that the Lie groups $G$ satisfying both conditions :
$G/K$ is Hermitian of tube type and $G/H$ is para-Hermitian, have a
nice description in terms of Euclidean Jordan algebras.
\bigskip

We shall briefly recall the link between Jordan algebras and the semi-simple
Lie groups considered above.

An algebra $V$ over $\mathbb R$ or $\mathbb C$ is said to be a
Jordan algebra if for all elements $x$ and $y$ in $V$, one has $
x\cdot y=y\cdot x$ and $ x\cdot(x^{2}\cdot y)=x^{2}\cdot(x\cdot y).$
For an element $x\in V$ let $L(x)$ be the linear map of $V$ defined
by $L(x)y:=x\cdot y$ and $P(x)=2L(x)^2-L(x^2)$ be the quadratic
representation of $V$. For $x$ and $y$ in $V$ one also defines an
endomorphism $D(x,y)$ of V given by $D(x,y)=L(xy)-[L(x),L(y)]$.

 We
denote by $\beta(x,y)$ the symmetric bilinear form on $V$ defined by
$\beta(x,y)={\rm Tr} L(x\cdot y)$.

Let $r$ and $n$ denote respectively the rank and the dimension of
the Jordan algebra $V$. The integer $d$ determined by $n=r+\frac
d2r(r-1)$ is called \emph{Peirce multiplicity}.
 For a regular element $x$, the minimal
polynomial $f_x$ is of degree $r$,
$$f_x(\lambda )=\lambda ^r -a_1(x)\lambda ^{r-1}+\cdots +(-1)^ra_r(x).$$
The coefficient $a_j$ is a homogeneous polynomial of degree $j$,
$\Delta (x):=a_r(x)$ is the \emph{Jordan determinant}, and ${\rm
tr}\, (x):=a_1(x)$ is the \emph{Jordan trace} of $x$.

A Jordan algebra $V$ is semi-simple if the form $\beta$ is
non-degenerate on $V$. A semi-simple Jordan algebra is unital, we
denote by $e$ its identity element.

A Jordan algebra $V_{o}$ over $\mathbb{R}$ is said to be Euclidean
if the bilinear form $\beta(x,y)$ is positive definite on $V_{o}$.

Let $V_{o}$ be an Euclidean Jordan algebra (EJA) from now one.  The
set
$$
\Omega:=\{x^{2}\:\vert\:x\:{\rm invertible}\:{\rm in}\: V_{o}\}
$$
is an open, convex, self-dual cone in $V_{o}$. Those properties of
$\Omega$ actually characterize $V_{o}$ as an EJA. The automorphism
group of $G(\Omega)$ of the cone $\Omega$ is defined by
$$
G(\Omega)=\{g\in GL(V_o)\,|\,g\Omega=\Omega\},
$$
and it is a reductive Lie group.

 Let $V$ be the complexification of
$V_{o}$. Consider the tube $T_{\Omega}=V_{o}+i\Omega\subset V$ and
the Lie group $G=Aut (T_{\Omega})$ of holomorphic automorphisms of
$T_{\Omega}$.
According to general theory \cite{[FK]} Ch. X. \S 5,
the group $G(\Omega)$ can be seen as a subgroup of $G$ as well as
the Jordan algebra $V_o$ it-self. Indeed, for every $u\in V_o$, the
translation $\tau_u:z\to z+u$ is a holomorphic automorphism of the
tube $T_\Omega$ and the group of all real translations $\tau_u$ is
an Abelian subgroup $N$ of $G$ isomorphic to the vector space $V_o$.

The subgroup of all affine linear transformations of the tube $P=
G(\Omega)\ltimes N$ is a maximal parabolic subgroup of $G$.

The subgroups $G(\Omega)$ and $N$ together with the inversion map
$j:x\To -x^{-1}$, generate the group $G$.

Let $\sigma$ be the involution of $G$ given by $ \sigma(g)=j\circ
g\circ j,\, g\in G.$ In the case when $V_o$ is a Euclidean Jordan
algebra this is a Cartan involution.  Let $K$ be a maximal compact
subgroup of $G$. Then the symmetric space $G/K\simeq T_{\Omega}$ is
an Hermitian symmetric space of tube type.

For $w\in V$ the endomorphism $D(w,\bar w)$ is Hermitian and one defines
an invariant spectral norm $|w|=\Vert D(w,\bar w)\Vert^{1/2}$. Let
$$
{\mathcal D}=\{w\in V\,:\,|w|<1\},
$$
be the open unit ball for the spectral norm. Then the Cayley transform $p:\,z\mapsto
(z-ie)(z+ie)^{-1}$ is a holomorphic isomorphism from the tube $T_\Omega$ onto the domain $\mathcal D$.
Thus the group of holomorphic automorphisms of $\mathcal D$ that one denotes $G({\mathcal D})=Aut({\mathcal D})$
is conjugate to $G$ : $G({\mathcal D})=pGp^{-1}$. We shall refer to the domain $\mathcal D$ as to the Harish-Chandra bounded realization of the symmetric space $G/K$.

We denote $\overline N=\sigma(N)$ and $\overline
P:=G(\Omega)\ltimes\overline N$.

From the geometric point of view the subgroup $\overline P$ can be
characterized in the following way:
$$
\overline P=\{g\in G'\:|\:g(0)=0\},
$$
where $G'$ is the subset of $G$ of all transformations well defined
at $0\in V_o$. It is open and dense in $G$. Moreover
$G'=NG(\Omega)\overline N$. The map $N\times
G(\Omega)\times\overline N\to G'$ is a diffeomorphism. We shall
refer to this decomposition as to the \emph{Gelfand-Naimark
decomposition} of the group $G$. Furthermore, for every
transformation $g\in G$ which is well defined at $x\in V_o$, the
transformation $gn_x$ belongs to $G'$ and its Gelfand-Naimark
decomposition is given by :
\begin{equation}\label{gelf-naimark}
gn_x=n_{g.x}(Dg)_x\bar n',
\end{equation}
where $(Dg)_x\in G(\Omega)$ is the differential of the conformal map
$x\to g.x$ at $x$ and $\bar n'\in\overline N$ (see \cite{[P]} Prop.
1.4).

\vskip 10pt

The flag variety ${\mathcal M}=G/\overline P$, which is compact, is
the \emph{conformal compactification} of $V_o$. In fact the map
$x\To(n_x\circ j)P$ gives rise to an embedding of $V_o$ into
$\mathcal M$ as an open dense subset, and every transformation in
$G$ extends to $\mathcal{M}$. \vskip 15pt

Let $\mathfrak g$ be the Lie algebra of the automorphism group $G$.
Euclidean Jordan algebras, corresponding Lie algebras of  infinitesimal automorphisms of tube domains, and their maximal compact
 subalgebras are given by the first table.

\begin{table}[ht]
 \begin{tabular}{|c|c|c|c|}
\hline
${\mathfrak g}$ & $\mathfrak k$ & $V$& $V_{o}$\\
\hline $\mathfrak{su}(n,n)$ &
$\mathfrak{su}(n)\oplus\mathfrak{su}(n)\oplus\Real$& $M(n,\mathbb
C)$ &
  $Herm(n,\mathbb C)$\\
$\mathfrak{sp}(n,\mathbb R)$& $\mathfrak{su}(n)\oplus\mathbb R$ &
$Sym(n,\mathbb C)$&
  $Sym(n,\mathbb R)$\\
$\mathfrak{so}^{*}(4n)$ & $\mathfrak{su}(2n)\oplus \mathbb R$ &
$Skew(2n,\mathbb C)$&
  $Herm(n,\mathbb H)$\\
$\mathfrak o(n,2)$ & $\mathfrak{o}(n)\oplus\mathbb R$ & $\mathbb
C^{n-1}\times\mathbb C$&
$\mathbb R^{n-1}\times\Real$\\
$\mathfrak{e}_{7(-25)}$ & $\mathfrak{e}_6\oplus\mathbb R$ &
$Herm(3,\mathbb O)\otimes\mathbb C$&
  $Herm(3,\mathbb O)$\\
\hline
\end{tabular}\\
\end{table}

Let us consider the involution $\eta$ of the complex Jordan algebra
$V$ given by $\eta (x+iy)=-x+iy\ (x,y\in V)$ and define the corresponding fix point sub-group
in $ G$ by $H =\{g\in G\, |\, \eta g\eta =g\}.$
 Clearly $H=G(\Omega)$. The involutions $\eta$
and $\sigma$ commute.

Notice that the involution we introduced is a particular case of a
conjugation of $V$ satisfying the following properties.
\begin{itemize}
\item $\eta\, (V_0)=V_0$,
\item $\eta\, (ie)=ie$,
\item $-\eta$ is a real Jordan algebra automorphism of $V$.
\end{itemize}

The factor space $G/H$ is a para-Hermitian symmetric space. It means
that its tangent bundle splits into two $G$-invariant sub-bundles
both isomorphic to the underlying Jordan algebra $V_o$.

We restricted all considerations to Euclidean Jordan algebras
therefore the para-Hermitian spaces $G/H$ that we get are of a
particular type, one calls them causal symmetric spaces of Cayley
type \cite{[FHO]}. Their infinitesimal classification is given in
the second table.
\begin{table}[ht]
\begin{tabular}{|c|c|c|c|}
  \hline
  $\mathfrak g$ &$ \mathfrak h$ &$ V$ &$ V_0$  \\
  \hline
$  \mathfrak{su}(n,n) $&$ \mathfrak{sl}(n,\mathbb{C})\oplus\mathbb{R}$ &$ {\rm M}(n,\mathbb{C})$ &$ {\rm Herm} (n,\mathbb{C}) $\\
$  \mathfrak{sp}(n,\mathbb{R})$ &$\mathfrak{sl}(n,\mathbb{R})\oplus\mathbb{R}$ &$ {\rm Sym}(n,\mathbb{C})$ &$ {\rm Sym }(n,\mathbb{R}) $\\
 $ \mathfrak{so}^*(4n)$&$\mathfrak{su}^*(2n)\oplus\mathbb{R} $&$ {\rm Skew}(2n,\mathbb{C})$ &$ {\rm Herm}(n,\mathbb{H}) $\\
$  \mathfrak{so}(n,2) $&$\mathfrak{so}(n-1,1)\times\mathbb R$ &$ \mathbb{C}^{n-1}\times\mathbb C $&$ \mathbb{R}^{n-1}\times\Real $\\
$  \mathfrak{e}_{7(-25)}$ &$\mathfrak{e}_{6(-26)}\oplus\mathbb R$ &$ {\rm Herm}(3,\mathbb{O})\otimes\mathbb{C} $&$ {\rm Herm}(3,\mathbb{O})$ \\
  \hline
\end{tabular}
\end{table}
\bigskip

\section{Two series of representations of $G$}
\subsection{Holomorphic discrete series}

Holomorphic induction from a maximal compact subgroup leads to a series of unitary representations of $G$,
called holomorphic discrete series representations, that one usually realizes on holomorphic sections of holomorphic vector bundles over $G/K$.

According to our convenience and easiness of presentation we shall use both bounded and unbounded
realizations of the symmetric space $G/K$. We start with the simplest case of scalar holomorphic discrete series.

 For a real parameter $\nu$ consider the weighted Bergman
spaces $H^{2}_{\nu}(T_{\Omega})$ of complex valued holomorphic
functions $f\in{\mathcal O}(T_{\Omega})$ such that
$$
\Vert f\Vert^{2}_{\nu}=\int_{T_\Omega}\vert
f(z)\vert^{2}\Delta^{\nu-2\frac nr}(y)dxdy<\infty,
$$
where $z=x+iy\in T_{\Omega}$. Note that the measure $\Delta^{-2\frac
nr}(y)dxdy$ on $T_{\Omega}$ is invariant under the action of the
group $G$. For $\nu>1+d(r-1)$
these spaces are non empty Hilbert
spaces with reproducing kernels. More precisely, the space
$H^{2}_{\nu}(T_{\Omega})$ has a reproducing kernel $K_\nu$ which is
given by
\begin{equation}\label{reprkernel}
K_\nu(z,w)=c_\nu\Delta\left(\frac{z-\bar w}{2i}\right)^{-\nu},
\end{equation}
where $c_\nu$ is some expression involving Gindikin's conical
$\Gamma$-functions (see \cite{[FK]} p.261).

 The action of $G$ on $H^{2}_{\nu}(T_{\Omega})$ given for every integer $\nu>1+d(r-1)$ by
\begin{equation}\label{hds}
\pi_{\nu}(g)f(z)={\rm Det}^{\nu}(D_{g^{-1}}(z))f(g^{-1}.z)
\end{equation}
is called a {\it scalar holomorphic discrete series representation}.\footnote
{Notice that in general one shows, by use of analytic continuation, that the reproducing kernel
(\ref{reprkernel}) is positive-definite
for a larger set of spectral parameters, namely for every $\nu$ in the so-called
Wallach set $W(T_\Omega)=\left\{0,\frac d2,\dots,(r-1)\frac d2\right\}\cup](r-1)\frac d2,\infty[$.
However we restrict our considerations only to the subset of $W(T_\Omega)$ consisting of integer
$\nu>1+d(r-1)$ in order to
deal with spaces of holomorphic functions.}

In the above formula $D_{g}(z)$ denotes the differential of the
conformal transformation $z\to g.z$ of the tube.\vskip 10pt

On the other hand side the corresponding action of the group $G(\mathcal D)$ can be realized as follows.
Let
$$
B(z,w)=1-D(z,w)+P(z)P(w),
$$
be the Bergman operator on $V$. Its determinant ${\rm det}B(z,w)$ is of the form $h(z,w)^{2n/r}$ where $h(z,w)$
is the so-called \emph{canonical polynomial}
(see \cite{[FK]} p.262). Notice that it is the pull back of $K_1(z,w)$ by the Cayley transform.

Then the group $G({\mathcal D})$ acts
on the space $H^2_\nu({\mathcal D})$ of holomorphic functions $f$ on $\mathcal D$ such that
$$
\Vert f\Vert_{\nu,\mathcal D}^2=c_{\nu}'\int_{\mathcal D}
|f(z)|^2h(z,z)^{\nu-2\frac nr}dxdy<\infty
$$
by the similar formula
$\pi_\nu(g)f(z)={\rm Det}^{\nu}(D_{g^{-1}}(z))f(g^{-1}.z).$\vskip 5pt

More generally let $\mathfrak g$ be the Lie algebra of the
automorphisms group $G({\mathcal D})$ with complexification ${\mathfrak g}_c$. Let
$\mathfrak g=\mathfrak k\oplus \mathfrak p$ be a Cartan
decomposition of $\mathfrak g$. Let $\mathfrak z$ be the center of
$\mathfrak k$. In our case the centralizer of $\mathfrak z$ in
$\mathfrak g$ is equal to $\mathfrak k$ and the center of $\mathfrak
k$ is one-dimensional. There is an element $Z_0\in\mathfrak z$ such
that $({\mbox {ad}}Z_0)^2=-1$ on $\mathfrak p$. Fixing $i$ a square
root of $-1$, one has ${\mathfrak p}_c={\mathfrak p}+i{\mathfrak
p}={\mathfrak p}_++{\mathfrak p}_-$ where ${\mbox {
ad}}Z_0\vert_{{\mathfrak p}_+}=i,\: {\mbox {ad}}Z_0\vert_{{\mathfrak
p}_-}=-i.$ Then
\begin{equation}
{\mathfrak g}_c={\mathfrak p}_+\oplus {\mathfrak k_c}\oplus
{\mathfrak p}_-.
\end{equation}
and $[{\mathfrak p}_{\pm},{\mathfrak p}_{\pm}]=0,[{\mathfrak
p}_+,{\mathfrak
 p}_-]={\mathfrak k}_c$ and $[{\mathfrak k}_c,{\mathfrak
p}_{\pm}]={\mathfrak p}_{\pm}.$ The vector space $\mathfrak p_+$ is isomorphic to $V$ and furthermore
it inherits its Jordan algebra structure. Let $G_c$ be a connected, simply
connected Lie group with Lie algebra ${\mathfrak g}_c$ and $K_c,
P_+$, $P_-, G,K,Z$ the analytic subgroups corresponding to ${\mathfrak
k}_c,{\mathfrak p}_+,{\mathfrak
 p}_-,{\mathfrak g}$, $\mathfrak k$ and $\mathfrak z$
respectively. Then $K_cP_-$ (and $K_cP_+$) is a maximal parabolic
subgroup of $G_c$ with split component $A=\exp i{\Bbb R}Z_0$. So the group $G=G(\mathcal D)_o$ is
closed in $G_c$.

Moreover, the exponential mapping is a diffeomorphism of ${\mathfrak
p}_-$ onto $P_-$ and of ${\mathfrak p}_+$ onto $P_+$
(\cite{[Helgas]} Ch.VIII, Lemma 7.8). Furthermore:
\begin{lem}\label{lem1}
a.The mapping $(q,k,p)\mapsto qkp$ is a diffeomorphism of $P_+\times
K_c\times P_-$ onto an open dense submanifold of $G_c$ containing
$G$.
\newline
b. The set $GK_cP_-$ is open in $P_+K_cP_-$ and $G\cap K_cP_-=K$.
\end{lem}
(see \cite{[Helgas]}, Ch VIII, Lemm\ae $\:$ 7.9 and 7.10).

Thus $G/K$ is mapped on an open, bounded domain $\mathcal D$ in
${\mathfrak p}_+$ This is an alternative description of the
Harish-Chandra bounded realization of $G/K$. The group $G$ acts on
$\mathcal D$ via holomorphic transformations.

Everywhere in this section we shall denote $\bar g$ the complex
conjugate of $g\in G_c$ with respect to $G$ (do not confuse with the
involution $\sigma$). Notice that $ P_+$ is conjugate to $P_-$.
\vskip .3cm For $g\in P_+K_cP_-$ we shall write
$g=(g)_+\:(g)_0\:(g)_-$, where $(g)_{\pm}\in P_{\pm},\:(g)_0\in
K_c$. For $g\in G_c,\: z\in\mathfrak p_+$ such that $g.\exp z\in
P_+K_cP_-$ we define
\begin{eqnarray}
\exp g(z)&=&(g.\exp z)_+\\
J(g,z)&=&(g. \exp z)_0.
\end{eqnarray}
$J(g,z)\in K_c$ is called the {\it canonical automorphic factor} of
$G_c$ (terminology of Satake).
\begin{lem}
\cite{[Satake]} Ch.II, Lemma 5.1. The map $J$ satisfies

(i) $J(g,o)=(g)_0$, for $g\in P_+K_cP_-$,

(ii) $J(k,z)=k$ for $k\in K_c, z\in\mathfrak p_+$.
\newline
If for $g_1,g_2\in G_c$ and $ z\in\mathfrak p_+$, $g_1(g_2(z))$ and
$g_2(z)$ are defined, then $(g_1g_2)(z)$ is also defined and

(iii) $J(g_1g_2,z)=J(g_1,g_2(z))\:J(g_2,z)$.
\end{lem}

 \vskip 10pt For
$z,w\in\mathfrak p_+$ satisfying $(\exp\bar w)^{-1}.\exp z\in
P_+K_cP_-$ we define
\begin{eqnarray}
K(z,w)&=&J((\exp\bar w)^{-1},z)^{-1}\\
&=&((\exp\bar w)^{-1}.\exp z)_0^{-1}.
\end{eqnarray}
This expression is always defined for $z,w\in\mathcal D$, for then
$$
(\exp\bar w)^{-1}.\exp
z\in\overline{(GK_cP_-)}^{-1}GK_cP_-=P_+K_cGK_cP_-=P_+K_cP_-.
$$
$K(z,w)$, defined on $\mathcal D\times\mathcal D$, is called the
{\it canonical kernel} on $\mathcal D$ ( by Satake). $K(z,w)$ is
holomorphic in $z$, anti-holomorphic in $w$, with values in $K_c$.
Here are a few properties:
\begin{lem}\label{lemK}
\cite{[Satake]}, Ch.II, Lemma 5.2. The map $K$ satisfies\\
(i) $K(z,w)=\overline{K(w,z)}^{-1}$ if $K(z,w)$ is defined,

(ii)$K(o,w)=K(z,o)=1$ for $z,w\in\mathfrak p_+$.
\newline
If $g(z),\bar g(w)$ and $K(z,w)$ are defined, then $K(g(z),\bar
g(w))$ is also defined and one has:

(iii) $K(g(z),\bar g(w))=J(g,z)\:K(z,w)\:\overline{J(\bar
g,w)}^{-1}$,
\end{lem}

\begin{lem}
\cite{[Satake]}, Ch.II, Lemma 5.3. For $g\in G_c$ the Jacobian of
the holomorphic mapping\newline $z\mapsto g(z)$, when it is defined,
is given by
$$
{\rm Jac}\:(z\mapsto g(z))={\mbox {Ad}}_{\mathfrak p_+}(J(g,z)).
$$
\end{lem}

For any holomorphic character $\chi:K_c\mapsto\Bbb C$ we define:
\begin{eqnarray}
j_{\chi}(g,z)&=&\chi(J(g,z)),\\
k_{\chi}(z,w)&=&\chi(K(z,w)).
\end{eqnarray}
Since $\chi(\bar k)=\overline{\chi(k)}^{-1}$ we have :
\begin{eqnarray}
k_{\chi}(z,w)&=&\overline{k_{\chi}(w,z)},\\
k_{\chi}(g(z),\bar
g(w))&=&j_{\chi}(g,z)k_{\chi}(z,w)\overline{j_{\chi}(\bar g,w)}
\end{eqnarray}
in place of Lemma (\ref{lemK}) (i) and (iii).

The character $\chi_1(k)=\det{\mbox {Ad}}_{\mathfrak p_+}(k),\:(k\in
K_c)$ is of particular importance. We call the corresponding
$j_{\chi_1},k_{\chi_1}$: $j_1$ and $k_1$. Notice that
\begin{equation}\label{jac}
j_1(g,z)=\det({\mbox {Jac}}\:(z\mapsto g(z))).
\end{equation}

Because of (\ref{jac}), $\vert k_1(z,z)\vert^{-1}d\mu(z)$, where
$d\mu(z)$ is the Euclidean measure on $\mathfrak p_+$, is a
$G-$invariant measure on $\mathcal D$. Indeed:
 \begin{eqnarray*}
 d\mu(g(z))&=&\vert j_1(g,z)\vert^2d\mu(z),\\
 k_1(g(z),g(z))&=&j_1(g,z)\:k_1(z,z)\:\overline{j_1(g,z)},\qquad{\rm
for}\: g\in G.
 \end{eqnarray*}
One can actually show that $k_1(z,z)>0$ on $\mathcal D$.
(\cite{[Satake]}, Ch.II, Lemma 5.8).

Let $\tau$ be an irreducible holomorphic representation of $K_c$ on
a finite dimensional complex vector space $W$ with scalar product
$\langle\:\vert\:\rangle$, such that $\tau_{\vert_K}$ is unitary.
\begin{lem}
For every $k\in K_c$ one has the identity $\tau^*(k)=\tau(\bar
k)^{-1}$.\end{lem}

This follows easily by writing $k=k_o\cdot\exp iX$ with $k_o\in K$,
$X\in \mathfrak k$ and using that $\tau_{\vert_K}$ is unitary.
\vskip 7pt
 Call $\pi_{\tau}=\Ind_K^G\tau$ and set
$W_{\tau}$ for the representation space of $\pi_{\tau}$. Then
$W_{\tau}$ consists of maps $f:G\mapsto W$ satisfying

(i) $f$ measurable,

(ii) $f(gk)=\tau^{-1}(k)f(g)$ for $g\in G,\:k\in K$,

(iii) $\int_{G/K}\Vert f(g)\Vert^2d\dot g<\infty$, \noindent where
$\Vert f(g)\Vert^2=\langle f(g)\vert f(g)\rangle$ and $d\dot g$ is
an invariant measure on $G/K$. Let us identify $G/K$ with $\mathcal
D$ and $d\dot g$ with $d_*z=k_1(z,z)^{-1}d\mu(z)$. Then $W_{\tau}$
can be identified with a space of maps on $\mathcal D$, setting
\begin{equation}\label{21}
\varphi(z)=\tau(J(g,o))f(g),
\end{equation}
if $z=g(o), f\in W_{\tau}$. Indeed, the right-hand side of
(\ref{21}) is clearly right $K-$invariant. The inner product becomes
$$
(\varphi\vert\psi)=\int_{\mathcal
D}\langle\tau^{-1}(J(g,o))\varphi(z)\vert\tau^{-1}(J(g,o))\psi(z)\rangle
d_*z.
$$
Since
$\tau^{-1}(J(g,o))^*\tau^{-1}(J(g,o))=\tau^{-1}(J(g,o)\overline{J(g,o)}^{-1})=\tau^{-1}(K(z,z))$
by Lemma (\ref{lemK}), we may also write
\begin{equation}\label{22}
(\varphi\vert\psi)=\int_{\mathcal
D}\langle\tau^{-1}(K(z,z))\varphi(z)\vert\psi(z)\rangle d_*z.
\end{equation}
The $G$-action on the new space is given by
\begin{equation}
\pi_{\tau}(g)\varphi(z)=\tau^{-1}(J(g^{-1},z))\varphi(g^{-1}(z)),\qquad(g\in
G,\:z\in\mathcal D).
\end{equation}
Now we restrict to the closed sub-space of holomorphic maps and call
the resulting Hilbert space ${\mathcal H}_{\tau}$. The space
${\mathcal H}_{\tau}$ is $\pi_{\tau}(G)$-invariant. We assume that
${\mathcal H}_{\tau}\not=\{0\}$.

The pair $(\pi_{\tau},{\mathcal H}_{\tau})$ is called a {\it vector-valued
holomorphic discrete series} of $G$.

In a similar way we can define the anti-holomorphic discrete series.
We therefore start with $\bar\tau$ instead of $\tau$ and take
anti-holomorphic maps. Then
\begin{equation}
\pi_{\bar\tau}(g)\psi(z)=\bar\tau^{-1}(J(g^{-1},z))\psi(g^{-1}(z)).
\end{equation}
for $\psi\in{\mathcal H}_{\bar\tau}$. One easily sees that
${\mathcal H}_{\bar\tau}=\bar{\mathcal H}_{\tau}$ and
$\pi_{\bar\tau}= \bar\pi_{\tau}$ in the usual sense. Notice that when the representation
$\tau$ is one dimensional we recover scalar holomorphic discrete series representations introduced above.

\vskip 10pt

 The Hilbert space ${\mathcal H}_{\tau}$ is
known to have a reproducing (or Bergman) kernel ${\mathcal
K}_{\tau}(z,w)$. Its definition is as follows. Set
$$
E_z:\varphi\mapsto\varphi(z)\qquad(\varphi\in {\mathcal H}_{\tau})
$$
for $z\in\mathcal D$. Then $E_z:{\mathcal H}_{\tau}\mapsto W$ is a
continuous linear operator, and ${\mathcal K}_{\tau}(z,w)=E_zE^*_w$,
being a ${\rm End}(W)$-valued kernel, holomorphic in $z$,
anti-holomorphic in $w$. In more detail :
\begin{equation}\label{25}
\langle\varphi(w)\vert\:\xi\rangle=\int_{\mathcal
D}\langle\tau^{-1}(K(z,z))\varphi(z)\vert\:{\mathcal
K}_{\tau}(z,w)\xi\rangle d_*z
\end{equation} for any $\varphi\in{\mathcal H}_{\tau},\:\xi\in W$ and
$w\in\mathcal D$.

Since ${\mathcal H}_{\tau}$ is a $G-$module, one easily gets the
following transformation property for ${\mathcal K}_{\tau}(z,w)$ :
\begin{equation}\label{26}
{\mathcal K}_{\tau}(g(z),g(w))=\tau(J(g,z)){\mathcal
K}_{\tau}(z,w)\tau(\overline{J(g,w)})^{-1}\qquad (g\in G,
z,w\in\mathcal D).
\end{equation}
Now consider $H(z,w)={\mathcal
K}_{\tau}(z,w)\cdot\tau^{-1}(K(z,w))$.

Clearly $H(g(z),g(w))=\tau(J(g,z))H(z,w)\tau^{-1}(J(g,z))$ for all
$z,w\in\mathcal D$. So, setting $z=w=o, g\in K$ we see that $H(o,o)$
is a scalar operator, and hence $H(z,z)=H(o,o)$ is so. But then
$H(z,w)=H(o,o)$. So, we get
\begin{equation}\label{27}
{\mathcal K}_{\tau}(z,w)=c\cdot\tau(K(z,w)),
\end{equation}
where $c$ is a scalar. The same reasoning yields that $\pi_{\tau}$
is {\it irreducible}. Indeed, if ${\mathcal H}\subset{\mathcal
H}_{\tau}$ is a closed invariant subspace, then $\mathcal H$ has a
reproducing kernel, say $K_{\mathcal H}$ and it follows that
$K_{\mathcal H}=c{\mathcal K}_{\tau}$, so either ${\mathcal
H}=\{0\}$ or $\mathcal H=\mathcal H_{\tau}$.\vskip10pt

Let us briefly recall the analytic realization of (some of)
vector-valued holomorphic discrete series representations of $G$. We
start with the irreducible representations of the maximal compact
subgroup $K$ which can be realized on the space of polynomials
${\mathcal P}(V)$ and which are parameterized by the weights
$\m=(m_1,\dots,m_r)\in\mathbb Z^r$ with $m_1\geq\dots\geq m_r\geq0$
and $m_1+\dots+m_r=m=|\m|$. These representations do not exhaust all
irreducible representations of $K$ but they will produce all
necessary components for our further discussion.

Let $V'$ be the dual vector space of $V\simeq\mathfrak p_+$. Consider the $m-th$ symmetric tensor power of $V'$.
It is naturally identified with the space ${\mathcal P}^m(V)$ of polynomials of degree $m$ on $V$. It is well known
(see for instance \cite{[FK90],[Schmid]}) that under the $K$-action this space decomposes multiplicity free into
a direct sum of irreducible sub-representations :
$$
{\mathcal P}^m(V)=\sum_{|\m|=m}^\oplus{\mathcal P}^{\m}(V),
$$
where ${\mathcal P}^{\m}(V)$ are irreducible representations of $K$ of highest weight $\m$.
This decomposition is often called the Kostant-Hua-Schmid formula and we refer the reader to the paper \cite{[FK90]}
for a precise description of spaces ${\mathcal P}^{\m}(V)$ and the corresponding highest weight vectors $\Delta_{\m}$.
We denote by $P_{\m}$ the orthogonal projection of ${\mathcal P}^m(V)$ onto ${\mathcal P}^{\m}(V)$.

Let $h(z,w)$ be as before the canonical polynomial on $V\times V$, then according to \cite{[FK90]}, for a real $\nu$
one has
$$
h^{-\nu}(z,w)=\sum_{\m}(\nu)_{\m}K_{\m}(z,w),
$$
where $K_{\m}(z,w)$ is the reproducing kernel of the space ${\mathcal P}^{\m}(V)$, and $(\nu)_{\m}$ stands for
the generalized Pochhammer symbol:
$$
(\nu)_{\m}=\prod_{j=1}^r\left(\nu-\frac d2(j-1)\right)_{m_j}=\prod_{j=1}^r\prod_{k=1}^{m_j}\left(\nu-\frac d2(j-1)+k-1\right).
$$
Denote ${\mathcal H}_\nu({\mathcal P}^{\m}(V))$ the Hilbert space of
holomorphic functions on $\mathcal D$ with values in ${\mathcal
P}^{\m}(V)$ admitting the reproducing kernel
$$
h^{-\nu}(z,w)\otimes^m K^t(z,w).
$$
Then, for an integer $\nu>1+d(r-1)$ and a given weight $\m$ the group $G$ acts on its unitarily and irreducibly by
\begin{equation}\label{pinum}
\pi_{\nu,\m}(g)f(z)=\Det(Dg^{-1}(z))^\nu\left(\otimes^m(dg^{-1})^t\right)\cdot f(g^{-1}.z),
\end{equation}
where $\otimes^m(dg^{-1})^t$ on ${\mathcal P}^{\m}(V)$ denotes the induced action of $(dg^{-1})^t$ on $V$.
\subsection{Maximal degenerate series}
Let $\Det(g)$ be the determinant of a linear transform $g\in
G(\Omega)\subset GL(V_o)$. We denote by $\chi(g)$ a particular
character of this reductive Lie group given by
$\chi(g):=\Det(g)^{\frac rn}$.

This character can be trivially extended to the whole parabolic
subgroup $\overline P$ by $\chi(h\bar n):=\chi(h)$ for every $h\in
G(\Omega),\:\bar n\in\overline N$.\vskip15pt

For every $\mu\in\mathbb C$ we define a character $\chi_{\mu}$ of
$\overline P$ by $ \chi_{\mu}(\bar p):=|\chi(\bar p)|^{\mu}. $

The induced representation $\pi_{\mu}^-=\Ind_{\overline
P}^{G}\left(\chi_{\mu}\right)$ of the group $G$ acts on the space
$$
\widetilde I_{\mu}:=\{f\in C^{\infty}(G)\:|\:f(g\bar
p)=\chi_{\mu}(\bar p)f(g),\forall g\in G,\bar p\in\overline P\},
$$
by left translations. A pre-Hilbert structure on $\widetilde
I_{\mu}$ is given by $\norm f^2=\int_K\abs{ f(k)}^2dk,$ where $K$ is
the maximal compact subgroup of $G$ associated with the Cartan
involution $\sigma$, and $dk$ is the normalized Haar measure of $K$.

According to the Gelfand-Naimark decomposition a function
$f\in\widetilde I_{\mu}$ is determined by its restriction
$f_{V_o}(x)=f(n_x)$ on $N\simeq V_o$. Let $ I_{\mu}$ be the subspace
of $C^{\infty}(V_o)$ of functions $f_{V_o}$ with $f\in \widetilde
I_{\mu}$. The group $G$ acts on $I_{\mu}$ by:
\begin{equation}\label{action}
     \pi_{\mu}^-(g)f(x)=|A(g,x)|^{\mu}f(g^{-1}.x),\:g\in G,\:x\in V_o,
\end{equation}
where $A(g,x):=\chi_{\mu}\bigl((Dg^{-1})_x\bigr).$ These
representations are usually called the \emph{maximal degenerate
series representations} of $G$.

One shows that the norm of a function $f(n_x)=f_{V_o}(x)\in I_{\mu}$
is given by:
\begin{equation}\label{norm}
     \norm f^2=\int_{V_o}\vert f_{V_o}(x)\vert^2h(x,-x)^{2\Re \mu+\frac{n}{r}}dx,
\end{equation}

where $h(z,w)$ is the \emph{canonical polynomial} introduced above.
Formula (\ref{norm}) implies that for $\Re \mu=-\frac{n}{2r}$ the
space $I_{\mu}$ is contained in $L^2(V_o)$ and the representation
$\pi_{\mu}^-$ extends as a unitary representation on $L^2(V_o)$.

Analogously the character $\chi$ can be extended to the subgroup
$P$ and one defines in a similar way the representation
$\pi^+_\mu=\Ind_{P}^{G}(\chi_{-\mu})$.

  Following the standard procedure
we introduce an intertwiner between $\pi_{\mu}^-$ and
$\pi_{\mu-\frac{n}{r}}^+$. Consider the map $\widetilde A _{\mu}$
defined on $\widetilde I_{\mu}$ by
\begin{equation}\label{intertwiner}
     f\To(\widetilde A_\mu f)(g):=\int_Nf(gn)dn,\quad\forall g\in G,
\end{equation}
where $dn$ is a left invariant Haar measure on $N$. One shows that
this integral converges for $\Re \mu>{n\over 2r_0}$.


\begin{prop}
For every $f\in\widetilde I_{\mu}$ the function $\widetilde A_\mu f$
belongs to $\widetilde I_{-\mu}$ and the map $\widetilde A_\mu$
given by (\ref{intertwiner}) intertwines the corresponding
representations of $G$:
\begin{equation}\label{intertw-repr}
     \widetilde{\pi}_{\mu-\frac{n}{r}}^+(g)(\widetilde A_\mu f)=\widetilde
A_\mu(\widetilde{\pi}_{\mu}^-(g)f),\:\forall
     f\in\widetilde I_{\mu},\:g\in G.
\end{equation}
\end{prop}

%
%

\section{Ring structures on the holomorphic discrete series}
In this section we discuss two different ring structures that one
can endow on the set of holomorphic discrete series.

\subsection{Laplace transform and the point-wise product}
We start with a generalization of a result on the usual point-wise
product due to A. and J. Unterberger (cf. \cite{[UU]} Lemma 3.1) in
the case when $G=SL(2,\Real)$.

\begin{thm} Let $V_0$ be a Euclidean Jordan
algebra and $T_\Omega$ be the corresponding tube domain
$V_0+i\Omega$. Consider two real numbers $\nu_1$ and $\nu_2$ such
that $\nu_1,\nu_2>1+d(r-1)=2\frac nr-1$ and two functions $F_1\in
H^2_{\nu_1}(T_\Omega)$ and $F_2\in H^2_{\nu_2}(T_\Omega)$. Then
their point-wise product $F_1\cdot F_2$ belongs to
$H^2_{\nu_1+\nu_2}(T_\Omega)$.
\end{thm}
In order to prove this statement recall the following result
(\cite{[FK]}, Theorem XIII 1.1 ). Let $\Gamma_\Omega$ denote the
Gindikin conical $\Gamma$-function.
\begin{lem}
Let $\nu$ be a real number, $\nu >{2n\over r} -1$. Let $L^2_\nu
(\Omega )$ be the space $L^2(\Omega ,\Delta (2u)^{-\nu + {n\over
r}}du)$. For any $f\in L^2_\nu (\Omega )$, set
\begin{equation}
F(z)=(2\pi )^{-n/2}\int_\Omega f(u) {\rm e}^{(z|u)}du.
\end{equation}
Then $F\in {H}^2_\nu (T_\Omega )$ and $f\mapsto F$ is a linear
isomorphism from $L^2_\nu (\Omega)$ onto ${ H}^2_\nu (T_\Omega )$.
Moreover
$$
\Vert F\Vert^2_\nu =\Gamma_\Omega \left(\nu - {n\over
r}\right)\,\Vert f\Vert^2_\nu.
$$
\end{lem}
Let now $F_1\in{ H}^2_{\nu_1}(T_\Omega ),\, F_2\in{
H}^2_{\nu_2}(T_\Omega )$ and let $u$ and $v$ correspond to $F_1$ and
$F_2$ respectively by the lemma, so $u\in L^2_{\nu_1}(\Omega ),\,
v\in L^2_{\nu_2}(\Omega )$. Then we shall show:
$$
\Vert u\ast v\Vert_{\nu_1 +\nu_2}\leq C(\nu_1 ,\nu_2 )\, \Vert
u\Vert_{\nu_1} \, \Vert v\Vert_{\nu_2},
$$
where $C(\nu_1,\nu_2 )$ is a constant. This is sufficient to prove
the theorem since the map $f\mapsto F$ sends convolutions to
point-wise products. Observe that $f$ is extended to $V$ by setting
it zero outside $\Omega$. More precisely we have

\begin{lem} Let $\nu_1 ,\nu_2>{2n\over r} -1,\ \Vert u\Vert_{\nu_1}<\infty, \, \Vert v\Vert_{\nu_2} <\infty$ for the measurable functions $u$ and $v$ on $\Omega$. Set
$$
(u\ast v)(\tau )=\int_{\Omega\cap (\tau -\Omega )} u(\tau -\eta )\,
v(\eta ) d\eta\quad (\tau\in\Omega ).
$$
Then $(u\ast v)(\tau )$ exists for almost all $\tau$, is measurable
and
$$
\Vert u \ast v\Vert_{\nu_1 +\nu_2}\leq C(\nu_1,\nu_2)\, \Vert
u\Vert_{\nu_1}\,\Vert v\Vert_{\nu_2}.
$$
\end{lem}
We only prove the estimate, since the rest of this lemma follows
from the same proof, applying Fubini's theorem at each step.
\par
We have to give an estimate for the integral
\begin{eqnarray*}
I&=&\int_\Omega (u\ast v )(\tau )\, \Delta (2\tau )^{{n\over r} -\nu_1 -\nu_2}\overline w(\tau )d\tau\\
&=& \int_\Omega\int_\Omega \Delta (2(\xi +\eta ))^{{n\over
r}-\nu_1-\nu_2} v(\eta) u(\xi )\overline w (\xi +\eta )\, d\xi d\eta
\end{eqnarray*}
under the assumption that $\Vert w\Vert_{\nu_1 +\nu_2}<\infty$
($w\in L^2_{\nu_1 +\nu_2}(\Omega )$).
\par
For any $t>0$, using the inequality
$$2 |u(\xi )v(\eta )|\leq t\, \Delta(2(\xi +\eta ))^{{1\over 2} (\nu_2 -\nu_1)}\, |u(\xi )|^2 +$$
$$t^{-1} \Delta (2(\xi +\eta ))^{{1\over 2} (\nu_1 -\nu_2)}\, |v(\eta )|^2$$
we get
\begin{eqnarray*}
2|I|&\leq & t\int_\Omega |u(\xi )|^2\, d\xi\int_{\xi +\Omega} \Delta (2\tau )^{{n\over r} -\nu_1 -\nu_2 +{1\over 2} (\nu_2 -\nu_1 )}\, |w(\tau )|\, d\tau\\
&+& t^{-1} \int_\Omega |v(\eta )|^2\, d\eta\int_{\eta +\Omega} \Delta (2\tau )^{{n\over r} -\nu_1 -\nu_2 +{1\over 2} (\nu_1 -\nu_2 )}\, |w(\tau )|\, d\tau\\
&\leq & \Vert w\Vert_{\nu_1 +\nu_2}\Big [ t \int_\Omega |u(\xi )|^2\, d\xi\, (\int_{\xi +\Omega} \Delta (2\tau )^{{n\over r} -2\nu_1}\ d\tau )^{1/2}\\
&+& t^{-1} \int_\Omega |v(\eta )|^2\, d\eta\, (\int_{\eta +\Omega}
\Delta (2\tau )^{{n\over r} -2\nu_2}\, d\tau )^{1/2}\Big ].
\end{eqnarray*}
Let us compute the expression $\int_{\xi +\Omega}\Delta (2\tau
)^{{n\over r} -2\nu_2} d\tau$ for $\xi\in\Omega$. Set $\xi =g\cdot
e$ for $g\in G(\Omega )$. The $G(\Omega )$-invariant measure on
$\Omega$ is equal to $\Delta (\tau )^{-{n\over r}} d\tau$, so that
we get
$$\int_{\xi +\Omega} \Delta (2\tau )^{{n\over r} - 2\nu_1} d\tau = 2^{n-2\nu_1 r}\int_{\xi +\Omega}\Delta (\tau )^{{2n\over r} -2\nu_1}\Delta (\tau)^{-{n\over r}}d\tau$$
\begin{equation}
=2^{n-2\nu_1 r}\int_{e +\Omega}\Delta (g\cdot \tau )^{{2n\over r}
-2\nu_1}\Delta (\tau)^{-{n\over r}}d\tau
\end{equation}
Now $\Delta (g\cdot \tau )=({\rm Det}\, g)^{r/n}\, \Delta (\tau
)=\Delta (g\cdot\tau )\,\Delta (\tau )$, so we obtain for the latter
expression
$$= \int_{e +\Omega} \Delta (2\tau )^{{n\over r} - 2\nu_1}\, d\tau\, \, \Delta (\xi )^{{2n\over r} - 2\nu_1}.$$
The term $\int_{e +\Omega} \Delta (2\tau )^{{n\over r} -2\nu_1}
d\tau$ has finally to be computed.
\par
We make the change of variable $\tau\mapsto \tau^{-1}$. Observe that
$(e +\Omega )^{-1}= (e -\Omega )\cap\Omega$. The differential of
$\tau\mapsto\tau^{-1}$ is $-P(\tau )^{-1}$ and $|{\rm Det}\,(-P(\tau
))^{-1}|=\Delta (\tau )^{2n\over r}$, see (\cite{[FK]},Prop. II 3.3
and Prop. III 4.2). So
$$\int_{e +\Omega} \Delta (2\tau )^{{n\over r} - 2\nu_1} d\tau = 2^{n - 2\nu_1 r} \int_{(e-\Omega )\cap\Omega} \Delta (\tau )^{-{n\over r} + 2\nu_1} \Delta (\tau )^{-{2n\over r}}d\tau$$
$$= 2^{n -2\nu_1 r}\int_{(e -\Omega )\cap\Omega}\Delta (\tau )^{-{2n\over r} + 2\nu_1}d\tau = 2^{n -2\nu_1 r} B_\Omega (-{2n\over r} + 2\nu_1 , {n\over r})$$
$$=2^{n - 2\nu_1 r} \displaystyle{\Gamma_\Omega (-{2n\over r} +2\nu_1 )\,\Gamma_\Omega ({n\over r})\over\Gamma_\Omega (-{n\over r} + 2\nu_1 )}.$$
So we obtain
$$|I| = \Big [ t \Vert u\Vert^2_{\nu_1} \, \Big\{ 2^{n - 2\nu_1 r} \displaystyle{\Gamma_\Omega (-{2n\over r} +2\nu_1 )\,\Gamma_\Omega ({n\over r})\over\Gamma_\Omega (-{n\over r} + 2\nu_1 )}\Big \}^{1/2} +$$
$$ t^{-1} \Vert v\Vert^2_{\nu_2}\, \Big \{2^{n -2\nu_2 r} \displaystyle{\Gamma_\Omega (-{2n\over r} +2\nu_2 )\,\Gamma_\Omega ({n\over r})\over\Gamma_\Omega (-{n\over r} + 2\nu_2 )}\Big\}^{1/2}\Big ] \, \Vert w\Vert_{\nu_1 +\nu_2}.$$
Taking the minimum for $t>0$, we get
$$\Vert u\ast v\Vert_{\nu_1 +\nu_2}\leq 2^{{n\over r} - (\nu_1 +\nu_2 ){r\over 2}} \Big \{\displaystyle{\Gamma_\Omega (-{2n\over r} +2\nu_1 )\,\Gamma_\Omega ({n\over r})\over\Gamma_\Omega (-{n\over r} + 2\nu_1 )}\Big\}^{1/4}\cdot$$
$$\Big\{\displaystyle{\Gamma_\Omega (-{2n\over r} +2\nu_2
)\,\Gamma_\Omega ({n\over r})\over\Gamma_\Omega (-{n\over r} +
2\nu_2 )}\Big\}^{1/4}\, \Vert u\Vert_{\nu_1}\,\Vert
v\Vert_{\nu_2}.$$

\begin{rem} If $G=SL(2,\Real)$, then $T_\Omega=\Pi$ and for every
$f\in H^2_\nu(\Pi)$ $\displaystyle\frac{df}{dz}\in
H^{2}_{\nu+2}(\Pi)$.
\end{rem}
\subsection{Product structure on $L^2(G/H)$}

There exists a $G$-equivariant embedding of square-integrable
functions on the causal symmetric space $G/H$
 into the composition algebra of Hilbert-Schmidt operators by means
of the following diagram:
$$
L^2(G/H)\hookrightarrow\pi_{\mu}^+\otimes\pi_{\mu}
^-\hookrightarrow\pi_{\mu}^+\otimes\overline{\pi_{\mu}^+} \simeq HS
(L^2(V_o),dx),$$ where $dx$ is the usual Lebesgue measure on $V_o$.
 \vskip 10pt

 The first arrow is of geometric nature and
it is given by the fact that the symmetric space $G/H$ is an open
dense subset of $G/P\cap G/\bar P$. The last isomorphism is given by
$$L^2(V_0,dx)\otimes\overline{L^2(V_0,dx)}\simeq
HS(L^2(V_0,dx).$$

This embedding  gives rise to a covariant symbolic calculus on
$G/H$.

In order to introduce the covariant symbolic calculus on $G/H$ we
start with the case of $G=SL(2,\Real )$ and
$H=\left\{\left(%
\begin{array}{cc}
  a & 0 \\
  0 & a^{-1} \\
\end{array}%
\right),\,a\in\Real^*\right\}$.

Let $P^-$ be the parabolic subgroup of $G$ consisting of the lower
triangular matrices
$$
P^-:\left(%
\begin{array}{cc}
  a & 0 \\
  c & a^{-1} \\
\end{array}%
\right),
$$
with $c\in\mathbb R, a\in\mathbb R^*$ and let $P^+$ be the group of
upper triangular matrices
$$
P^+:\left(%
\begin{array}{cc}
  a & b \\
  0 & a^{-1} \\
\end{array}%
\right),
$$
with $b\in\mathbb R,a\in\mathbb R^*$. The group $G$ acts on the
sphere $\displaystyle S=\left\{ s\in\mathbb R^2:\,\Vert
s\Vert^2=1\right\}$ and acts transitively on $\widetilde{S}=S/\sim$,
where $s\sim s'$ if and only if $s=\pm s'$, by
$$
g.s=\frac{g(s)}{\Vert g(s)\Vert}.
$$
 Clearly ${\rm Stab}(\widetilde{0,1})=P^-$. So $\widetilde
S\simeq G/P^-$. Similarly $\widetilde S\simeq G/P^+$: $\widetilde
S=G.(\widetilde{1,0})$. If $ds$ is the usual normalized surface
measure on $S$, then
$$
d(g.s)=\Vert g(s)\Vert^{-2}ds.
$$
For $\mu\in\mathbb C$, define the character $\omega_\mu$ of $P^\pm$
by
$$
\omega_\mu(p)=|a|^\mu.
$$
Consider $\displaystyle \pi_\mu^{\pm}={\rm
Ind}_{P^\pm}^G\omega_{\mp\mu}$.

Both $\pi_\mu^+$ and $\pi_\mu^-$ can be realized on
$C^{\infty}(\widetilde S)$, the space of smooth functions $\phi$ on
$S$ satisfying
$$
\phi(-s)=\phi(s),\quad (s\in S).
$$
The formula for $\pi_\mu^-$ is
$$
\pi_\mu^-(g)\phi(s)=\phi(g^{-1}.s)\Vert g^{-1}(s)\Vert^{\mu}.
$$
Let $\theta$ be the Cartan involution of $G$ given by
$\theta(g)={}^tg^{-1}$. Then
$$
\pi_\mu^+(g)\phi(s)=\phi(\theta (g^{-1}).s)\Vert
\theta(g^{-1})(s)\Vert^{\mu}.
$$

Since here $$\theta\left(%
\begin{array}{cc}
  a & b \\
  c & d \\
\end{array}%
\right)=w\left(\begin{array}{cc}
  a & b \\
  c & d \\
\end{array}%
\right)w^{-1}
$$
with $w=\left(\begin{array}{cc}
  0 & 1 \\
  -1 & 0 \\
\end{array}%
\right),$ one has that $\pi_\mu^-\sim\pi_\mu^+$.

Let $(\,,\,)$ denote the standard inner product on $L^2(S)$:
$$
(\phi,\psi)=\int_S\phi(s)\overline{\psi(s)}ds.
$$
Then this form is invariant with respect to the pairs
$$
(\pi_\mu^-,\pi_{-\bar\mu-2}^-),\quad{\rm and}\quad
(\pi_\mu^+,\pi_{-\bar\mu-2}^+).
$$
Therefore if $\Re\mu=-1$, then the representations $\pi_\mu^\pm$ are
unitary, the inner product being $(\,,\,)$.

$G$ acts also on $\widetilde{S}\times\widetilde{S}$ by
\begin{equation}\label{1}
g.(u,v)=(g.u,\theta(g)v).
\end{equation}
This action is \underline{not} transitive: the orbit
$$(\widetilde{S}\times\widetilde{S})^\#_s=\{(u,v):\,\langle
u,v\rangle\neq0\}=G.((\widetilde{0,1}),(\widetilde{0,1}))
$$
is dense. Moreover $(\widetilde{S}\times\widetilde{S})^\#_s\simeq
G/H$.

The map
$$
f\rightarrow f(u,v)|\langle u,v\rangle|^{-1+is},\quad (s\in\mathbb
R),
$$
is a unitary $G$-isomorphism between $L^2(G/H)$ and
$$\displaystyle\pi_{-1+is}^-\hat\otimes_{{}_2}\pi_{-1+is}^+$$ acting on
$L^2(\widetilde{S}\times\widetilde{S})$. The latter space is
provided with the usual inner product.

Define the operator $A_\mu$ on $C^\infty(\widetilde{S})$ by the
formula
$$
A_\mu\phi(s)=\int_S|\langle s,t\rangle|^{-\mu-2}\phi(t)dt.
$$
This integral is absolutely convergent for $\Re\mu<-1$, and can be
analytically extended to the whole complex plane as a meromorphic
function. It is easily checked that $A_\mu$ is an intertwining
operator
$$
A_\mu\pi_\mu^{\pm}(g)=\pi_{-\mu-2}^{\mp}(g)A_\mu.
$$
The operator $A_{-\mu-2}\circ A_\mu$ intertwines $\pi_\mu^\pm$ with
itself, and is therefore a scalar $c(\mu)$. It can be computed using
$K$-types:
$$
c(\mu)=\pi\frac{\Gamma\left(\frac{\mu+1}2\right)\Gamma\left(-\frac{\mu+1}2\right)}{\Gamma\left(\frac{-\mu}2\right)
\Gamma\left(1+\frac{\mu}2\right)}.
$$
One also shows that : $A_\mu^*=A_{\bar\mu}$. So that for $\mu=-1+is$
we get (by abuse of notation):
$$
c(s)=\pi\frac{\Gamma\left(\frac{is}2\right)\Gamma\left(-\frac{is}2\right)}{\Gamma\left(\frac{1-is}2\right)\Gamma\left(\frac{1+is}2\right)},
$$
and moreover
$$
A_{(-1+is)}\circ A_{(-1+is)}^*=c(s)I,
$$
so that
$\pi^{-\frac12}\frac{\Gamma\left(\frac{1+is}2\right)}{\Gamma\left(\frac{is}2\right)}
A_{(-1+is)}=d(s)A_{(-1+is)}$ is a unitary intertwiner between
$\pi^+_{-1+is}$ and $\pi^-_{-1-is}$.

We thus get a $\pi^-_{-1+is}\hat\otimes_{{}_2}\bar\pi^-_{-1+is}$
invariant map from $L^2(G/H)$ onto
$L^2(\widetilde{S}\times\widetilde{S})$ given by
\begin{eqnarray*}
f&\rightarrow&d(s)\int_Sf(u,w)|\langle
u,w\rangle|^{-1+is}|\langle v,w\rangle|^{-1-is}dw\\
&=&(T_sf)(u,v),\qquad s\neq 0.
\end{eqnarray*}
This integral does not converge: it has to be considered as obtained
by analytic continuation.

Define the product $\,f\#_s g\,$ for $f,g\in L^2(G/H)$ as follows:
$(T_sf)(u,v)$ is the kernel of a Hilbert-Schmidt operator $Op(f)$ on
$L^2(\widetilde S)$. Then we set:
$$
Op(f\#_s g)=Op(f)\circ Op(g).
$$
We then have:

$\bullet\, \Vert f\sharp_s\,g\Vert_2\leq\Vert f\Vert_2\cdot\Vert g\Vert_2.$

$\bullet\, Op(L_x f)=\pi^-_{-1+is}(x)Op(f)\pi^-_{-1+is}(x^{-1})$, so
$$
L_x(f\#_s g)=(L_xf)\#_s (L_xg),
$$
for $x\in G$.

Let us write down a formula for $f\#_s g$; we have:
\begin{eqnarray}\label{sharp}
&&d^{-2}(s)(f\#_s
g)(u,v)\\&=&\int_S\int_Sf(u,x)g(y,v)|[u,y,x,v]|^{-1+is}d\mu(x,y),
\end{eqnarray}
where $d\mu(x,y)=|\langle x,y\rangle|^{-2}dxdy$ is a $G$-invariant
measure on $\widetilde S\times\widetilde S$ for the $G$-action
(\ref{1}). Here
$$
[u,y,x,v]=\frac{\langle u,x\rangle\langle y,v\rangle}{\langle
u,v\rangle\langle y,x\rangle}.
$$

For a generic causal symmetric space of Cayley type $G/H$ the
composition formula of two symbols $f,g\in L^2(G/H)$ is defined in a
similar way. In order to keep a reasonable size of this note we just
indicate the flavor of the explicit formula. Recall that $G/H$ is
para-Hermitian, $T_{eH}(G/H)\simeq V_0\oplus V_0$, and hence
functions on it can be seen as functions on $V_0\times V_0$.
Therefore $f\sharp_s\,g$ is given by a double integral on $V_0$, of
$f$ and $g$ against an appropriate power of the quotient of four
functions that are integral kernels of the intertwining operator
(\ref{intertwiner}), exactly as in (\ref{sharp}).

\subsection{Product structure on $L^2(G/H)_{\rm hol}$}

One says that a symmetric space $G/H$ has discrete series
representations if the set of representations of $G$ on minimal closed
invariant subspaces of $L^2(G/H)$ is nonempty. According to a
fundamental result of Flensted-Jensen \cite{[FJ]} the discrete
series for $G/H$ is nonempty and infinite if
$$rank(G/H)=rank(K/K\cap H).$$

For a causal symmetric space of Cayley type $G/H$ this condition is
fulfilled and on can realize part of its discrete series as
holomorphic discrete series representations of the group $G$.

More precisely assume that $\pi$ is a scalar holomorphic discrete series
representation of $G$, i.e. it acts on ${\mathcal
H}_\pi\subset\mathcal O(D)\cap L^2(D,dm_\pi)$ where $D$ is some symmetric domain (the
image of the tube $T_\Omega$ by the Cayley transform) and where $dm_\pi(w)$ is a measure on $D$ associated to $\pi$.
In such a
case the Hilbert space ${\mathcal H}_\pi$ has a reproducing kernel
$K_\pi(z,w)$.

Assume that the representation $\pi$ occurs as a multiplicity free
closed subspace in the Plancherel formula for $L^2(G/H)$ (actually
this is the case in our setting, see \cite{[OO]} Theorem 5.9).

Consider $\xi_{\pi}\in{\mathcal H}_\pi^{-\infty}$ the unique up to
scalars $H$-fixed distribution vector associated to $\pi$ (see
\cite{[OO]} p. 142 for the definition of
$\xi_{\pi}=\phi_{\lambda}(z)$). It gives rise to a continuous
embedding map
$$
{\mathcal J}_\pi:{\mathcal H}_\pi\hookrightarrow
L^2(G/H)\subset{\mathcal D}'(G/H)
$$
given for any analytic vector $v\in{\mathcal H}_\pi^{\infty}$ by
\begin{equation}\label{j}
({\mathcal J}_\pi v)(x)=\langle v,\pi(x)\xi_{\pi}\rangle,\qquad x\in
G/H,
\end{equation}
where by abusing notations we write $\pi(x)$ instead of $\pi(g)$
with $x=g.H\in G/H$.

For any fixed $w\in D$ let us define the function
$v_w:=K_\pi(\cdot,w)$ which is actually a real analytic vector in
${\mathcal H}_\pi$.

Consider now the following function :
$$
g_w(x):=({\mathcal J}v_w)(x),\qquad x\in G/H,\:w\in D.
$$
\vskip 10pt

Because of the reproducing property of the Hilbert space ${\mathcal
H}_\pi$ for every $f\in{\mathcal H}_\pi$ one can write
$$
f(z)=\int_DK_\pi(z,w)f(w)dm_\pi(w).
$$

Furthermore, if such a function $f$ is an analytic vector for the
representation $\pi$ : $f\in {\mathcal H}_\pi^{\infty}$, then
$$
({\mathcal J}_\pi f)(x)=\int_Df(w)g_w(x)dm_\pi(w).
$$

Choosing an appropriate normalization in (\ref{j}) one can get the
embedding ${\mathcal J}_\pi$ isometric. Therefore the subspace
generated by $g_w(x)$,$\:w\in D$ is a closed subspace of $L^2(G/H)$
isometric to some holomorphic discrete series representation of $G$.
( see \cite{[OO]} Theorem 5.4 for the precise statement).

The dual map ${\mathcal J}_\pi^*:{\mathcal D}(G/H)\mapsto{\mathcal
H}_\pi$ is defined by
\begin{eqnarray*}
\langle {\mathcal J}_\pi^*\phi,f\rangle&=&\langle\phi,{\mathcal
J}_\pi
f\rangle\\&=&\int_{G/H}\int_D\phi(x)\overline{f(w)g_w(x)}dm(w)d\nu(x),\:\forall\phi\in{\mathcal
D}(G/H),
\end{eqnarray*}
where $d\nu(x)$ denotes the invariant measure on $G/H$. Therefore we
have,
$$
({\mathcal J}_\pi^*\phi)(w)=\int_{G/H}\phi(x)
\overline{g_w(x)}d\nu(x).
$$
Similar observations are valid for vector-valued holomorphic discrete series
representations as well.

Define the set
$$
L^2(G/H)_{\rm hol}=\bigoplus_{\pi\in\hat G'_{\rm hol}}{\mathcal
J}_\pi({\mathcal H}_\pi)
$$
where $\hat G'_{\rm hol}$ denotes the set of equivalence classes of
unitary irreducible holomorphic discrete series representations of
$G$ with corresponding character $\tau$ trivial on $H\cap Z$.
Notice that the space $L^2(G/H)_{hol}$ decomposes multiplicity
free into irreducible subspaces (\cite{[OM]}).

\bigskip
 According to
\cite{[OO]} the $H$-fixed distribution vector $\xi_k=\xi_{\pi_\nu}$,
associated with the scalar holomorphic discrete series
representation $\pi_\nu$ (see (\ref{hds})) is given up to a constant
by
\begin{equation}\label{xi}
\xi_k(z)=\Delta\left(\frac{\eta(z)-\bar z}{2i}\right)^{-\frac
\nu2},\quad z\in V_0+i\Omega.
\end{equation}
\bigskip

{\bf Example.} If $G=SL(2,\Real)$ the holomorphic discrete series of $G$
are only scalar and according to (\ref{hds}) act on $H^2_k(\Pi)\,
(k\geq 2, k\in\mathbb N)$ by
$$
\pi_k(g)f(z)=(cz+d)^{-k}f\left(\frac{az+b}{cz+d}\right),\qquad
g^{-1}=\left(%
\begin{array}{cc}
  a & b \\
  c & d \\
\end{array}%
\right).
$$

The involution $\eta$ is given here by $\eta(z)=-\bar z$, the subgroup $H$ is isomorphic to $SO(1,1)$ and
according to (\ref{xi}) $\xi_k(z)=(-\bar z)^{-\frac k2}$ what
corresponds
 precisely to $\frac k2$-th power of the Unterbergers generating function (\cite{[UU]} Prop 3.3.)
 $$
 g_z(s,t)=\frac{s-t}{(s-\bar z)(t-\bar z)},
 $$
 evaluated at the base point $(0,\infty)$ of the orbit $G/H$.

In this case it is well known \cite{[OO]} that only $\pi_k$ with $k$
{\bf{ even}} can be uniquely realized on $L^2(G/H)$. Therefore we have
$$
L^2(G/H)_{hol}=\bigoplus_{k\,{\rm even}}{\mathcal J}_k(H^2_k(\Pi)),
$$

In \cite{[UU]} Theorem 3.6 the authors showed that the set $
L^2(G/H)_{hol}$ is closed under the non commutative product $\#_s$
(\ref{sharp}). They give an explicit formula for the components of
$f\#_s g$ for $f\in {\mathcal J}_k(H^2_k(\Pi))$ and $g\in{\mathcal
J}_\ell(H^2_\ell(\Pi))$ in terms of Rankin-Cohen brackets
(\ref{rcb}). The method they developed is elegant but technical and
from our point of view not well adapted for generalization.
\bigskip

To get more insight in the product structure of $L^2(G/H)_{\rm
hol}$, we rely on some recent results by T.Kobayashi
(\cite{[Kobayashi]},Theorem 7.4),

We are going to show that $L^2(G/H)_{hol}$ is closed under the
product $\#_s$. It is, because of the continuity of the product,
sufficient to show the following theorem.

\begin{thm}\label{main}
Let ${\mathcal H}_\pi$ and ${\mathcal H}_{\pi'}$ be two irreducible
closed subspaces of $L^2(G/H)_{\rm hol}$. Then
$$
{\mathcal J}_\pi(f)\#_s {\mathcal J}_{\pi'}(g)\in L^2(G/H)_{\rm hol}.
$$
for every $f\in \mathcal{H}_\pi$ and $g\in\mathcal{H}_{\pi'}$.
\end{thm}

Even in the case of $G=SL(2,\mathbb R)$ this result reduces the
computations of \cite{[UU]} in an interesting way. The proof of this
theorem follows from a fairly recent result by Kobayashi saying:

\begin{thm}\label{koba}(Kobayashi). Let $\pi$ and $\pi'$ be holomorphic discrete
series representations of $G$. Then the representation
$$
\pi\hat\otimes_{{}_2}\pi'
$$
 decomposes discretely into holomorphic discrete series representations of $G$ with finite multiplicities. Moreover,
$\pi\hat\otimes_{{}_2}\pi' $ is $K$-admissible, i.e. every
irreducible representation of $K$ occurs in it with finite
multiplicity.
\end{thm}

In general we do not have a multiplicity free decomposition.

Now let us show our theorem. The map \\ $f\otimes g\to {\mathcal J}_\pi(f)\#_s {\mathcal J}_{\pi'}(g)$
clearly gives rise to a $K-$ and
$\mathcal{U}(\mathfrak{g})-$equivariant linear map
$$
(\mathcal{H}_\pi\hat\otimes_{{}_2}\mathcal{H}_{\pi'})^K=\mathcal{H}_\pi^K\otimes
\mathcal{H}_{\pi'}^K \rightarrow L^2(G/H),
$$ and thus the result
follows for $f$ and $g$ $K-$finite, and then, by continuity of the
product, for all $f$ and $g$.\bigskip

\textbf{Example} The decomposition of the tensor product of two holomorphic discrete
series for $SL(2,\mathbb R)$ was obtained by J.Repka
(\cite{[Repka]}) in full generality using the Harish-Chandra-modules
techniques, and it is given by
$$
\pi_n\hat\otimes_{{}_2}\pi_m=\bigoplus_{k=0}^\infty\pi_{m+n+2k},
$$
so that this reduces the computations in \cite{[UU]} even more.

At the same time V.F.Molchanov obtained the same result
(decomposition of all possible tensor products of unitary
irreducible representations of $SO(2,1)$) in \cite{[Molchanov]}. He
realized such tensor products on functions defined on the
one-sheeted hyperboloid and gave all the Fourier coefficients of the
positive-definite kernel that defines the Hilbert structure of these
unitary representations.
\vskip 15pt
In the general situation we have to consider also
vector-valued holomorphic discrete series representations.
Indeed, according to the theorem (\ref{koba}) and particularly to the result stated in the theorem 3.3
in \cite{[PZ]} the tensor product of two scalar holomorphic discrete series representations decomposes multiplicity free in the direct sum of unitary irreducible vector-valued holomorphic discrete series representations:
$$
{\mathcal H}_{\nu_1}\otimes{\mathcal H}_{\nu_2}=\sum_{\m\geq0}{\mathcal H}_{\nu_1+\nu_2}({\mathcal P}^{\m}(V')),
$$
in the case when $\nu_1\geq\nu_2>1+ d(r-1)$.

In order to understand the previous decomposition we have to identify its different ingredients.

First, we see an element of the tensor product $
{\mathcal H}_{\nu_1}\otimes{\mathcal H}_{\nu_2}$ as a holomorphic function $F(z,w)$ on $D\times D$. Therefore, according
to \cite{[Range]} Coroll. 6.26, p. 269 for any positive integer $m$ one can write a Taylor expansion formula:
$$
F(z,w)=\sum_{j=0}^m(F^{(j)}(z),\otimes^j(w-z))+(F^{(m+1)}(z,w),\otimes^{m+1}(z-w)),
$$
where $F^{(j)}(z)$ are ${\mathcal P}^j(V')$-valued holomorphic functions on $D$, $F^{(m+1)}(z,w)$
is a ${\mathcal P}^{m+1}(V')$-valued holomorphic function on $D\times D$ uniquely determined by the data of $F(z,w)$, and
$(\,,\,)$ denotes the standard pairing of corresponding vector spaces.

Second, according to Peetre \cite{[Peetre]} consider an $End(V)$-valued holomorphic differential form on $D$ defined for
every fixed $w_1,w_2\in V$ and $z\in D$ by
$$
\Omega(z;,w_1,w_2)=d_zB(z,w_1)B(z,w_1)^{-1}-d_zB(z,w_2)B(z,w_2)^{-1},
$$
and denote by $\omega(z;w_1,w_2)$ it trace $-\frac r{2n}{\mathrm{tr}}\Omega(z;w_1,w_2)$. The former differential form plays a crucial role in the construction of intertwining operators for tensor products.

Namely,
for fixed $w_1$ and $w_2$ the expression $$h(z,w_1)^{-\nu_1}h(z,w_2)^{-\nu_2}P_{\m}\otimes^{|\m|}\omega(z;w_1,w_2)$$ can be seen as an element
of the space $\overline{{\mathcal H}_{\nu_1}}\otimes\overline{{\mathcal H}_{\nu_2}}$ dual of ${\mathcal H}_{\nu_1}\otimes{\mathcal H}_{\nu_2}$. Let
$\langle\,,\,\rangle$ stand for the corresponding pairing.
Then the operator $I_{\m}$ given by
\begin{equation}\label{im}
I_{\m}(f\otimes g)(z)=\langle h(z,\cdot)^{-\nu_1}h(z,\cdot)^{-\nu_2}P_{\m}\otimes^{|\m|}\omega(z;\cdot,\cdot),
f\otimes g\rangle,
\end{equation}
is a $G$-equivariant map from $(\pi_{\nu_1}\otimes\pi_{\nu_2},{\mathcal H}_{\nu_1}\otimes{\mathcal H}_{\nu_2})$ to the space of ${\mathcal P}^{\m}
(V)$-valued holomorphic functions on $D$ seen as the representation space of $\pi_{\nu_1+\nu_2,\m}$ (see (\ref{pinum}).

Theorem 4.4 in \cite{[PZ]} which is an extended version
of the main theorem in \cite{[Peetre]} gives a description of this map.
Summarizing and using Theorem (\ref{main}), we get
\begin{prop} Let $\nu_1\geq\nu_2>1+d(r-1)$ and $f\in{\mathcal H}_{\nu_1},\,g\in {\mathcal H}_{\nu_2}$.
Assume that $\left({\mathcal H}_{\nu_1}^{-\infty}\right)^H$ and $\left({\mathcal H}_{\nu_2}^{-\infty}\right)^H$ are not
reduced to $\{0\}$.
Then
$$
{\mathcal J}_{\pi_{\nu_1}}(f)\#_s\,{\mathcal J}_{\pi_{\nu_2}}(g)=\sum_{\m\geq 0}c_{\m,s}
{\mathcal J}_{\pi_{\nu_1+\nu_2}}(B_{\m,\nu_1,\nu_2}(f,g)),
$$
where $c_{\m,s}$ are fundamental constants given by the $\#_s$ product of the reproducing kernels of the corresponding
Bergman spaces ${\mathcal H}_{\nu_1}$ and ${\mathcal H}_{\nu_2}$ and $B_{\m,\nu_1,\nu_2}$ is such a bi-differential operator on ${\mathcal H}_{\nu_1}\otimes {\mathcal H}_{\nu_2}$ that
\begin{equation}\label{product}
I_{\m}(B_{\m,\nu_1,\nu_2}(f,g))=
\sum_{|\n|+|\n'|=m}C_{|\m|}^{|\n|}\cdot\frac{(-1)^{|\n|}}{(\nu_1)_{\n}(\nu_2)_{\n'}}
\cdot
P_{\m}\left(P_{\n}\partial^{|\n|}f\otimes P_{\n'}\partial^{|\n'|}g\right),
\end{equation}
with $\n$ and $\n'$ being all possible
weights such that $|\n|+|\n'|=|\m|$.
\end{prop}

In case when the group under consideration is  $SU(1,1)$ formula (\ref{product}) reduces to the Rankin-Cohen brackets initially introduced for the $SL(2,\Real)$-action. The fact that the expression remains the same
in both
compact and non-compact realizations of the Riemannian symmetric space $G/K$, is due to the fact that the groups
$SU(1,1)$ and $SL(2,\Real)$ are real forms of the same complex Lie group $SL(2,\mathbb C)$ and therefore covariant
differential operators on $ SU(1,1)/S(U(1)\times U(1))$ and $SL(2,\Real)/SO(2,\Real)$ are isomorphic via the analytic
continuation in the complexification of these symmetric spaces. This phenomenon holds in general for covariant differential operators on $D$ and on $T_\Omega$ and for transvectants in particular as was noticed by Peetre \cite{[Peetre]} p. 1076.

For this reason it would be natural to call the bi-differential operators occurring in (\ref{product}) {\emph{generalized Rankin-Cohen brackets.}}

{\bf Open questions}
\begin{itemize}
\item The construction we described is valid for holomorphic discrete series representations with
spectral parameter $\nu>1+d(r-1)$. However it would be interesting to understand whether this can be extended to the whole Wallach set.
\item A possible relationship with Vertex algebras, already mentioned in \cite{[Peetre]} and \cite{[Zagier1]}, were pointed out
to us by I.Cherednik. It is a challenge to investigate this link.
\item We have seen that the ring structure on
$L^2(G/H)_{\rm hol}$ is related to the tensor product of holomorphic
discrete series representations. Does it reflect, via the
Tannaka-Krein duality, the existence of a certain Hopf algebra that
would govern the non-commutative product $\sharp_s$ ?
\item According to the  Beilinson-Bernstein classification of
$(\mathfrak g,K)$-modules one can study
representations of semi-simple Lie groups in
terms of $D$-modules on associated flag varieties.
Can the ring structure on the set of
holomorphic discrete series be
interpreted as a
cup-product on the sheaves associated with closed $H$-orbits on the flag variety?
\end{itemize}
\bibliographystyle{amsplain}

\end{document}